\documentclass[11pt,leqno]{article}
\usepackage{latexsym,enumerate}
\usepackage{amsmath,amssymb}
\usepackage{amsthm}
\title{On the Sylow graph of a group and Sylow normalizers\footnote{The second and third authors have been supported by
Proyecto MTM2007-68010-C03-03, Ministerio de Educaci\'{o}n y Ciencia and FEDER, Spain.\newline The first author thanks the Universitat de Val\`{e}ncia and the Universidad Polit\'{e}cnica de Valencia for their
warm hospitality during the preparation of this paper.}}
\author{L. S. Kazarin\\
\footnotesize{Department of Mathematics, Yaroslavl P. Demidov State University}\\
 \footnotesize{Sovetskaya Str 14, 150014 Yaroslavl, Russia}\\
 \footnotesize{E-mail: Kazarin@uniyar.ac.ru}\\
 \\
A. Mart\'{\i}nez-Pastor\\
\footnotesize{Escuela T\'{e}cnica Superior de Ingenier\'{\i}a Inform\'{a}tica,}\\
 \footnotesize{Instituto Universitario de Matem\'{a}tica Pura y Aplicada IUMPA-UPV}\\
 \footnotesize{Universidad Polit\'{e}cnica de valencia}\\
\footnotesize{Camino de Vera, s/n,  46022 Valencia, Spain}\\
 \footnotesize{E-mail: anamarti@mat.upv.es}\\
 \\
 and \\
 \\
 M. D. P\'{e}rez-Ramos\\
\footnotesize{Departament d'\`{A}lgebra, Universitat de Val\`{e}ncia}\\
 \footnotesize{C/ Doctor Moliner 50, 46100 Burjassot
({Val\`{e}ncia}), Spain}\\
\footnotesize{E-mail: Dolores.Perez@uv.es}}
\date{}
\newtheorem{thm}{Theorem}[section]
\newtheorem{lem}[thm]{Lemma}
\newtheorem{prop}[thm]{Proposition}
\newtheorem*{teor}{Main Theorem}
\newtheorem{cor}[thm]{Corollary}
\theoremstyle{definition}
\newtheorem*{de}{Definition}
\newtheorem{ex}{Example}

\theoremstyle{remark}
\newtheorem{rem}[thm]{Remark}
\newtheorem*{rem*}{Remark}
\DeclareMathOperator{\n}{n}
\DeclareMathOperator{\Char}{Char}
\DeclareMathOperator{\GF}{GF}
\DeclareMathOperator{\modul}{mod}
\DeclareMathOperator{\Syl}{Syl}
\DeclareMathOperator{\Aut}{Aut}
\DeclareMathOperator{\Autdiag}{Autdiag}

\DeclareMathOperator{\Inndiag}{Inndiag}
\DeclareMathOperator{\Sz}{Sz}
\DeclareMathOperator{\di}{d}
\DeclareMathOperator{\E}{E}

\begin{document}
\maketitle

\begin{abstract}
Let $G$ be a finite group and $G_p$ be a Sylow $p$-subgroup of $G$
for a prime $p$ in $\pi(G)$, the set of all prime divisors of the
order of $G$. The automiser  $A_p(G)$ is defined to be the group
$N_G(G_p)/G_pC_G(G_p)$. We define the Sylow graph $\Gamma_A(G)$ of
the group $G$, with set of vertices $\pi(G)$, as follows: Two
vertices $p,q\in\pi(G)$ form an edge of $\Gamma_A(G)$ if either
$q\in\pi(A_p(G))$ or $p\in \pi(A_q(G))$. The following result is
obtained:

Theorem: Let  $G$ be a finite almost simple group. Then the graph $\Gamma_A(G)$ is connected
and has diameter at most $5$.

We also show how this result can be applied to derive information
on the structure of a group from the normalizers of its Sylow
subgroups.
\end{abstract}

\section{Introduction}
All groups considered are finite.

Let $G$ be a finite group and $G_p$ be a Sylow $p$-subgroup of $G$
for any $p\in \pi(G)$, where $\pi(G)$ denotes the set of prime
divisors of the order of $G$. Then the automiser  $A_p(G)$ is
defined to be the group $N_G(G_p)/G_pC_G(G_p)$. We construct the
Sylow graph $\Gamma_A(G)$ of the group $G$, with set of vertices
$\pi(G)$, as follows: Two vertices $p,q\in\pi(G)$ form an edge of
$\Gamma_A(G)$ if either $q\in\pi(A_p(G))$ or $p\in \pi(A_q(G))$.
It is convenient to write $p\rightarrow q$ if $q\in \pi(A_p(G))$.
For almost simple groups, that is, for groups with a unique
minimal normal subgroup, which is a non-abelian simple group, the
following result is obtained in the paper.
\begin{teor} Let  $G$ be an almost simple group. Then the graph $\Gamma_A(G)$ is connected and has diameter at most $5$.
\end{teor}
It follows from this theorem that often (with few exceptions) the
nor\-malizer of a Sylow $p$-subgroup $G_p$ of a finite almost
simple group $G$ is larger than $G_pC_G(G_p)$, for $p$ an odd
prime.  In particular this is always true if $G$ is a finite
simple group (see \cite{TV, Gur}). However, we obtain here more
precise information in graph terms for almost simple groups.

The study of the Sylow  graph of a group was first motivated in
relation with the influence of the Sylow normalizers (normalizers
of the Sylow sub\-groups) on the structure of finite groups. In
Section~\ref{Sylnorm} we show how the result on the Sylow graph of
a group can be applied to derive some information about the group
within this framework. Examples 1, 2 and 3 show certain classes of
groups, which can be considered as extensions of nilpotent groups,
and which have the property that a group belongs to the class if
and only if its Sylow normalizers belong to the class
(Corollaries~\ref{cor0} and \ref{cor2-3}). For the classes in
Example 1 the problem has been studied  in \cite{DDG08} under the
assumption of the result on the Sylow graph of a group as a
conjecture; though it appears to be a gap in the proof.

Our results extend some theorems of Glauberman~\cite{G}, Bianchi
et al.~\cite{H},  and D'Aniello et al.~\cite{DGP}. The
generalization in \cite{H} of a result in ~\cite{G}  states that a
group is nilpotent whenever its Sylow normalizers are nilpotent.
In \cite{DGP} a similar result is
obtained for classes of groups with nilpotent Hall subgroups for
given sets of primes, in the universe of all finite soluble groups.
In fact, within this universe, this result provides a
characterization of the subgroup-closed saturated formations whose
elements are characterized by the Sylow normalizers belonging to
the class. The corresponding classes are called
covering-formations of soluble groups (see Definition in
Section~\ref{Sylnorm}).

Covering-formations form a special family of the so-called
nilpotent-like Fitting formations. These are non-empty saturated Fitting
formations locally defined by $\mathcal E_{\pi(p)}$, the class of
all finite $\pi(p)$-groups for a set of primes $\pi(p)$ containing $p$, for every prime $p$ in the characteristic of the formation. Different restrictions on the
sets of primes $\pi(p)$ provide extensions of the nilpotent groups
from different points of view. We refer to Ballester-Bolinches et
al. \cite{BMPP} for an account of this development. The classes of
groups in Examples 1, 2 and 3 are  particular types of
covering-formations.
\section{Preliminaries}
Here we collect some lemmas which will be useful in the sequel.

Let $n$ be a positive integer and $p$ a prime number. A prime $r$ is said
to be \textit{primitive with respect to the pair $\{p, n\}$} if
$r$ divides $p^n-1$ but $r$ does not divide $p^e-1$ for every
integer $e$ such that $1\le e< n$. It was proved by
Zsigmondy~\cite{Z} that such a primitive prime $r$ exists  unless either $n=2$ and $p$ is a Mersenne prime or
$\{p,n\}=\{2,6\}$. Moreover, it is well-known that for this prime $r$, it holds
 $r-1\equiv 0 \pmod n$ and, in particular, $r\ge n+1$.

\begin{lem} \emph{(Frattini argument)}
\label{L0} Let $G$ be a group, $N$ be a normal subgroup of $G$ and
$S$ be a Sylow subgroup of $N$. Then $G=NN_G(S)$.
\end{lem}

For notation  in this section we refer to \cite{Ca72}. Let $G$ be
a connected simple algebraic group with  root system $\Theta
=\Theta (G)$ over an algebraically closed field of
characteristic $r>0$. For every simple group $L$ of Lie type
defined over the field $\GF(r^t)$, split or twisted, there exists $G=G_{ad}$
and a Frobenius map $\sigma $ of $G$ such that
$L=O^{r'}(G_\sigma)$, where $G_\sigma=\{x\in G\mid \sigma(x)=x\}$.
The inner-diagonal automorphism group $\hat L$ of $L$ coincides
with $G_\sigma$. If $G$ has a root system $\Theta $, we  say that
$L$ is of type $\Theta $ and write $L=\Theta (r^t)$ (see also
\cite{Ca85}).

\begin{lem}\label{TV}  \emph{\cite[Theorem 3.3]{TV}}
Let $L=O^{r'}(G_\sigma)$ have a root system $\Theta \not=A_l\
(l>1),\ D_{2l+1},\ E_6$. Then every semisimple element of $L$ of
odd prime order is conjugate to its inverse, i.e.,  it is real.
\end{lem}

In view of Lemma~\ref{TV} we will call a group of Lie type
\emph{fine} if it has no root system in the set $\{A_l\ (l>1),\
D_{2l+1},\ E_6\}$.
\begin{lem}\label{L2.4}  A field automorphism of a finite simple group of Lie type centralizes
the Weyl subgroup of the group modulo its Cartan subgroup.
\end{lem}
\noindent\textbf{Proof.} We use some facts from \cite{Stein}.  Let $L$ be a finite simple group of Lie type $\Theta $ over the field $F:=\GF(q)$.  $L$
has a generating system of elements of
the form $x_\alpha(\lambda )$ with $\alpha\in \Theta,~\lambda \in F$. The
elements $w_\alpha(\lambda)=x_\alpha(\lambda)x_{-\alpha}(-\lambda^{-1})x_\alpha(\lambda)$
generate the monomial subgroup $N$ and the elements of the form
$h_\alpha(\lambda)=w_\alpha(\lambda)w_\alpha(1)^{-1}$ generate the Cartan
subgroup $H$ of $L$, which is normal in $N$. Moreover there exists
a coset representatives of $H$ in $N$ consisting of elements
$\omega_\alpha=w_\alpha(1)$. If $\phi$ is an automorphism of the
field $F$, then the action of the
extension of $\phi$ on $L$ can be described as
$x_\alpha(\lambda)^\phi=x_\alpha(\lambda^\phi)$. It is easy to see that
$\omega_\alpha^\phi=\omega_\alpha$, which proves the result.\qed

\begin{lem}\label{L*}  Let $L$ be a simple group of Lie type over a field of
characteristic $r>0$ and let $P$ be a Sylow $p$-subgroup of $L$, for a prime $p\not=r$.
Then:

(a) $P$ normalizes some maximal torus of $L$.

(b)  If $p$ does not divide the order of the Weyl group of $L$, then $P$ is
contained in some maximal torus of $L$ and hence $P$ is abelian.
\end{lem}
\noindent\textbf{Proof.} See  \cite[Corollary II.5.19]{BCCI}. \qed

\begin{lem}\label{L**} \emph{(\cite{Bour}, \cite{Ca72})} Let  $L$ be  a simple group of Lie type and let $W$ be its Weyl group. Then:

If $L$ is of type $A_n$,
then  $W\cong S_{n+1}$ (of order $(n+1)!$).

If $L$ is  of  type $D_n$, then $W$ is a semidirect product of an elementary abelian $2$-group of order
$2^n$ and a symmetric group $S_n$.

If $L$ is  of  type $^2D_n$, then $W$ is a semidirect product of an elementary abelian $2$-group of order
$2^{n-1}$ and a symmetric group $S_{n-1}$.

If $L$ is  of type $E_6$, then $W$
has a subgroup of index $2$ isomorphic  to $P\Omega_6^-(2)\cong
PSp_4(3)$.

If $L$ is  of  type $^2E_6$, then $W$ is
isomorphic to the Weyl group of $F_4$ of order $2^73^2$.
\end{lem}
\section{Properties of the graph $\Gamma (G)=\Gamma_A (G)$}

We write $\Gamma (G)$ instead of $\Gamma_A (G)$ for brevity, as only this type of graphs are studied in the paper.
The set of edges of $\Gamma (G)$ is denoted by $\E(\Gamma(G))$.
For a group $G$ and a prime $p$, $\Syl_p(G)$ denotes the set of all Sylow
$p$-subgroups of the group $G$.
\begin{lem}\label{L3} Let $G_0$ be a normal subgroup of a group $G$ such that
$|G/G_0|=p\in \pi(G)$ and let $r$ be a prime, $r\not=p$. If $r\rightarrow s\in \E(\Gamma(G_0))$, then  $r\rightarrow s\in \E(\Gamma(G))$.
\end{lem}
\noindent
\textbf{Proof.} Let $R$ be a Sylow $r$-subgroup of $G_0$. Then
$|A_r(G_0)|\equiv 0\,(\modul s)$. It follows from Lemma~\ref{L0} that
$|N_G(R)|=|N_{G_0}(R)|p$ and  $s\in \pi(A_r(G))$.\qed

\begin{lem}\label{L4b} Let $G$ be a group and $P\in \Syl_p(G)$, $p>2$. If there exists an element $1\not=z\in Z(P)$
which is real in $G$, then
 $p\rightarrow 2\in \E(\Gamma(G))$.
\end{lem}
\noindent \textbf{Proof.} The extended centralizer of
$z$ in $G$ is defined to be
$$H:=\{y\in G\mid z^y= z^{\pm 1}\}.$$
Clearly $H$ is a subgroup of $G$ and $|H:C_G(z)|=2$. Since $P\leq C_G(z)< H$ it is obvious that $|A_p(H)|\equiv 0\,(\modul 2)$.
Moreover $\Syl_p(H)\subseteq \Syl_p(G)$ implies that $\pi(A_p(H))\subseteq \pi(A_p(G))$ and then the result follows.\qed

\bigskip
For a group $G$, define the \emph{distance} $\rho (p,q)$, for
primes $p,q\in \pi(G)$, as the minimal number $k$ for which there
is a chain $p_0=p, p_1,\ldots,p_k=q$ such that every pair $(p_i,
p_{i+1})$  forms an edge. The largest $\rho (p,q)$, $p,q\in
\pi(G)$, $p\not=q$, is the {\it  diameter} of $\Gamma (G)$,
denoted by $\di(\Gamma (G))$.

\begin{cor}\label{C2} The graph $\Gamma(S_n)$ of a symmetric group
of degree $n>2$ is connected and $\di(\Gamma (S_n))\le 2$.
\end{cor}

\noindent \textbf{Proof.} This follows by Lemma~\ref{L4b} since
each element of a symmetric group is real.\qed

\begin{lem}\label{L3.7} Let $G$ be an almost simple group and $N$ be  a non-trivial normal
subgroup of $G$ such that $(|G/N|,|N|)=1$. Then $$\di(\Gamma
(G))\leq \di(\Gamma (N))+2.$$

If for each pair of different primes $r, s\in \pi(G/N)$ there
exists $p\in \pi(N)$ such that $rs$ divides $|A_p(G)|$, then
$$\di(\Gamma (G))\leq \di(\Gamma (N))+1.$$
\end{lem}
\noindent \textbf{Proof.} Suppose that $r\in \pi(G)\setminus
\pi(N)$.  By Lemma~\ref{L0} we have that $|G|_r$ divides
$|N_G(G_p)|$ for any prime $p\in \pi(N)$ and $G_p\in \Syl_p(G)$.  Since $G$ is almost simple, it is clear that $|G|_r$ does not divide
$|C_G(G_p)|$ for some prime $p\in \pi(N)$, and then $p\rightarrow r\in \E(\Gamma(G)$. Hence, if $r,s\in \pi(G)\setminus
\pi(N)$, it follows that $\rho (r,s)\leq \di(\Gamma (N))+2$.

 If $rs$ divides $|A_p(G)|$ for a prime $p\in \pi(N)$, then $\di(\Gamma
(G))\leq \di(\Gamma (N))+1$, since $\rho (r,p)=1$ and $\rho (s,p)=1$.\qed

\section{Proof of the Main Theorem}
Corollary~\ref{C2} is a first step towards the proof of our Main
Theorem, which states that the Sylow graph of an almost simple
group is connected and has diameter at most $5$. In the next
results we consider separately different cases, depending on the
simple group appearing as the socle of an almost simple group, and
prove that in all cases the desired result holds.
\subsection*{Alternating and sporadic groups}
\begin{thm}\label{P1} Let $G=A_n$ be an alternating group of degree $n\ge 5$.
Then the graph $\Gamma(G)$ is connected and $\di(\Gamma(G))\le 3$.
If $n\not =p,\, p+1$, with $p$ a prime, $p\equiv 3\,(\modul 4)$, then
$\di(\Gamma(G))= 2$.
\end{thm}
\noindent \textbf{Proof.} Let $p$ be a prime,  $3\le p\le n$.  Assume that $p\rightarrow 2\not\in
\E(\Gamma (G))$. By Lemma~\ref{L4b} there exists an element  $z$ of order $p$ which is not real in $G$. If $p\leq n-2$, then
$z=a_1a_2\ldots a_k$ where $a_1,a_2,\ldots ,a_k$ are disjoint
cycles of length $p$, $k\geq 1$ and $kp+r=n$.

Suppose first that $r\geq 2$. In this case there exists an odd
permutation $\sigma $ centralizing $z$. On the other hand,
$z^g=z^{-1}$ for some $g\in S_n$. If $g\in S_n\setminus A_n$, then
$z^{\sigma g}=z^{-1}$ and $\sigma g\in A_n$. But this is a
contradiction by Lemma~\ref{L4b}. Hence $r\leq 1$.

If $k>1$, then there exists an involution  $\sigma \in
S_n\setminus A_n$ permuting $a_1$ and $a_2$ and fixing all letters
not appearing in $a_1$ and $a_2$. This gives a contradiction as
above  since $z^\sigma =z$ and $z^{\sigma g}=z^{-1}$.

Assume that either $n=p$ or $n=p+1$. There is an element $g\in S_n$
inverting $z$ which is a product of $\frac{p-1}{2}$ transpositions
$g=(1,p)(2,p-1)\ldots (\frac{p-1}{2},\frac{p+3}{2})$, as a
permutation of the set $\{1,2,\ldots ,n\}$. Since $g\in
S_n\setminus A_n$, we have $\frac{p-1}{2}\equiv 1\,(\modul 2)$,
which implies $p\equiv 3\,(\modul 4)$.

In this case there exists an odd prime $s$ dividing $p-1$. It is
easy to see that $p\rightarrow s\in \E(\Gamma (G))$ and $s\rightarrow
2\in \E(\Gamma (G))$ by  using arguments as above. Therefore for each pair
$p,r\in \pi(G)$ there is a path from $p$ to $r$ of length at most
$3$. The theorem is proved.\qed

\begin{rem*} The  alternating group $A_8$ shows that the estimate of
Theorem~\ref{P1} is sharp. In this case $\pi(A_8)=\{2,3,5,7\}$ and
$\E(\Gamma(A_8))=\{7\rightarrow 3,  3\rightarrow 2, 5\rightarrow 2 \}$.
\end{rem*}

\begin{thm}\label{S} Let $L$ be a sporadic simple group and $L\leq G\leq \Aut(L)$.
Then the graph $\Gamma(G)$ is connected and $\di(\Gamma(G))\le 5$.
\end{thm}
\noindent \textbf{Proof.} This is easy to prove using the
description in \cite{GL} of the normalizers of Sylow subgroups of
the considered groups.

Observe that there are only $3$ sporadic simple groups whose Sylow
graphs have the diameter exactly $5$. They are the Mathieu group
$M_{23}$, the Baby Monster group $B$ and the Monster group $M$.\qed
\subsection*{Almost simple groups with fine socle}
We recall that  a group of Lie type is called \emph{fine} if it
has no root system in the set $\{A_l\ (l>1),\ D_{2l+1},\ E_6\}$.

\begin{thm}\label{fs} Let $G$ be an almost simple group whose socle is a fine group of Lie type.
Then the graph $\Gamma (G)$ is connected and $\di(\Gamma(G))\le 4$.
\end{thm}
\noindent \textbf{Proof.} Let  $L$ be a fine group of Lie type
over the field $\GF(r^t)$ and set
$G_0=G\cap \Inndiag(L)$.
\smallskip

Note that, by Lemma~\ref{TV}, for a group $X$ such that $L\leq
X\leq G$  every semisimple element of odd prime order in $L$ is
real. Suppose that $p$ is an odd prime in $\pi(L)\setminus \{r\}$.
Then a Sylow $p$-subgroup $L_p$ of $L$ is normal in some Sylow
$p$-subgroup $X_p$ of $X$ and there exists an element $1\not= z\in
Z(X_p)\cap L_p $ such that $z$ is conjugate with $z^{-1}$. By
Lemma~\ref{L4b} we have that $p\rightarrow 2\in \E(\Gamma(X))$.

\smallskip
Assume  first  that the characteristic $r$ of $L$ is odd.

\smallskip
Let $L\not=L_2(q)$ with $q\equiv 3\pmod 4$. By the above property
and inspecting the orders of the Borel subgroups of $L$,  we have
that $p\rightarrow 2\in \E(\Gamma(G_0))$ for any odd prime $p\in
\pi(G_0)$. In particular $\di(\Gamma(G_0))\leq 2$.

This fact also holds for any group $G_1\geq G_0$ such that
$\pi(G_1)=\pi(G_0)=\pi(L)$. Hence we consider $G_1\lhd G$ such
that $\pi(G/G_1)\cap \pi(L)=\emptyset$.

We claim that $\di(\Gamma(G))\le 3$. Now we
may apply Lemma~\ref{L3.7}. Since the graph automorphism of $L$ does not
centralize any Sylow $r$-subgroup of $G_1$, we have that $\di(\Gamma
(G))\leq \di(\Gamma (G_1))+1\leq 3$.
\smallskip

Let now $L=L_2(q)$ where $q\equiv 3\pmod 4$. It is well-known that
the normalizer  of a Sylow $r$-subgroup of $L$  is a
Frobenius group of order $\frac{q(q-1)}{2}$, and the centralizer of
any field automorphism of $L$ does not contain a Sylow $r$-subgroup
of $L$. If $|G/L|\equiv 0\pmod 2$, then again $r\rightarrow 2\in \E(\Gamma(G))$.  Clearly, if $v\in \pi (t)\cap (\pi (G)\setminus
\{r\})$, then also $r\rightarrow v$. It is obvious  that $\di(\Gamma(G))\le 3$ if $PGL_2(q)$
is a subgroup of $G$,  and $\di(\Gamma(G))\le 4$ if $|G/L|$ is odd.

The case  $L=L_2(q)$ with $q=3^{3^a}$ demonstrates the following
interesting phenomena. If $G$ is an extension of $L$ by the group
of order dividing $3^a>1$, then the Sylow $3$-subgroups of $G$ are
self-normalizing (see, for instance, \cite{Gur}). If $K$ is the
normalizer in $G$ of a Sylow $2$-subgroup of $G$ (and of $L$),
then $K/O(K)\cong A_4$, the alternating group of degree $4$, and
we have $2\rightarrow 3$. In this case  $\di(\Gamma(G))=2$. In the
case when $2$ divides the order of $G/L$, we have  also
$\di(\Gamma(G))=2$, but the structure of the graph $\Gamma(G)$ is
slightly different.
\smallskip

Now consider the case when $r=2$ and $t>1$.

Again,  it follows as above that $p\rightarrow 2\in
\E(\Gamma(G_1))$, for each odd prime $p\in \pi(L)$,  and a group
$G_1$ such that $G_0\leq G_1\unlhd G$, $\pi(G_1)=\pi(G_0)=\pi(L)$
and $\pi(G/G_1)\cap \pi(L)=\emptyset$. In particular, $\di(\Gamma
(G_1))\leq 2$.

Suppose that $p\in \pi(G)\setminus \pi(L)$. This is the case only
when $p$ is the order of some field automorphism $\varphi $ of
$L$. Since $\varphi $  centralizes no Sylow $2$-subgroup of $L$, it
follows from Lemma~\ref{L3.7} that $\di(\Gamma (G))\leq \di(\Gamma
(G_1))+1\leq 3$, and the theorem is proved.\qed
\subsection*{Almost simple groups with socle $L=E_6^\pm (q)$}
In  Lemmas~\ref{LE1}~--~\ref{LE4} and Theorem~\ref{LE5}, $L$ is
$E_6^\pm (q)$, $r$ denotes the characteristic of $L$ and $W(L)$
its Weyl group. Here $E_6^+(q)$ and $E_6^-(q)$ stand for
$E_6(q)$ and $^2E_6(q)$, respectively. It is convenient to write
$L=E_6^\epsilon =E_6^\pm , \epsilon =\pm $.

If $d=(3,q-\epsilon 1)$, then the order of $L$ is as follows:
\begin{align*}|L|=|E_6^\pm  (q)|= &\frac{1}{d}q^{36}
(q-\epsilon 1)^6(q+\epsilon 1)^4(q^2-\epsilon q+1)^2(q^2+\epsilon q+1)^3\\
&(q^2+1)^2(q^4-q^2+1)(q^4+1)\frac{q^5-\epsilon 1}{q-\epsilon 1}
(q^6+\epsilon q^3+1).
\end{align*}

It was proved by Stensholt  in \cite{St1, St2} that $F_4(q)$ is a
subgroup of $L$. The order of the Weyl group of $F_4(q)$ is
$2^73^2$, and all semisimple elements of odd prime order in
$F_4(q)$ are real by Lemma~\ref{TV}, since this is a fine group.

\begin{lem}\label{LE1}
If $s\geq 5$ is a prime  dividing $q^4-q^2+1$, $q^2\pm q+1$,
$q^2+1$ or $q^4+1$, then $s\rightarrow 2\in \E(\Gamma(L))$.
\end{lem}
\noindent \textbf{Proof.}  If $s>5$ the Sylow $s$-subgroups of $L$
are abelian by Lemma~\ref{L*}.

Consider first the case $L=E_6^+(q)$. Since all semisimple
elements of odd prime order in $F_4(q)$ are real, we have that
$s\rightarrow 2\in \E(\Gamma (L))$ by Lemma~\ref{L4b}.

Suppose that $1\not= R\in \Syl_5(L)$. Then $5$ divides  $q^2+1$.
In this case $R$ normalizes a torus $T$ such that
$|N_L(T)/T|\equiv 0\,(\modul 5)$. It follows from \cite{CaLN} that
there is just one such torus $D_4(a_1)$ and its index in the
normalizer is $96$. Hence  the Sylow $5$-subgroups  are abelian and $5\rightarrow 2\in \E(\Gamma(L))$.

Let now $L=E_6^-(q)$. Since the Weyl group of $L$ is of type $F_4$
of order $2^73^2$, it follows that every Sylow $s$-subgroup  is
abelian, and since $s\in\pi(F_4(q))$, every element of order $s$ in $L$ is
real. Hence $s\rightarrow 2\in \E(\Gamma(L))$ as asserted. \qed

\begin{lem}\label{LE2} If $s\not=5$ is a prime dividing $q^4+\epsilon q^3+q^2+\epsilon q+1$ and not dividing
$q-\epsilon 1$,  then
$s\rightarrow 5\in \E(\Gamma(L))$.
\end{lem}
\noindent {\bf Proof.} It is easy to see that $s\mid q^5-\epsilon
1$ and all other factors $q^i\pm 1$ dividing $|L|$ have common
prime divisors with $q^5-\epsilon 1$ only in $\pi(q-\epsilon 1)$.
By Lemma~\ref{L*} the Sylow $s$-subgroups of $L$ are
abelian.

Consider first the case $L=E_6^+(q)$. By \cite{CaLN} there are
tori $T$ of types $A_4$ and $A_4+A_1$ of  $L$, of orders
$\frac{(q^5-1)(q-1)}{(q-1,3)}$ and $\frac{(q^5-1)(q+1)}{(q-1,3)}$,
respectively. In both cases $|N_L(T):T|=10$.

Note that $L$ contains a section $K$ isomorphic to $L_2(q)\times
L_5(q)$, which follows from the Dynkin diagram of this group. $K$ has a Frobenius subgroup of order  $\frac{q^5-1}{(q-1)(q-1,5)}\cdot 5$  and the Sylow $s$-subgroups of $L$ are cyclic.
Hence there is an edge $s\rightarrow 5\in \E(\Gamma(L)$.

In the case $L= E_6^-(q)$, the group $L$ contains a
section $K$ isomorphic to $L_2(q)\times U_5(q)$ and $U_5(q)$ has a
Frobenius subgroup of order $5(q^4-q^3+q^2-q+1)$. Again  the Sylow
$s$-subgroups of $L$ are cyclic. Hence there is an edge
$s\rightarrow 5\in \E(\Gamma(L)$, which concludes the proof.\qed
\begin{lem}\label{LE3} $5\rightarrow 2\in \E(\Gamma(L))$; if $r\not=2$, then $ r\rightarrow 2\in \E(\Gamma(L))$.
\end{lem}
\noindent {\bf Proof.} It follows from the description of the
Borel and Cartan subgroups of $L$ that $r\rightarrow 2\in
\E(\Gamma(L))$ if $r\not=2$. Hence we may assume that $5\not=r$
and $5$ divides either $q^2-1$ or $q^2+1$. By Lemma~\ref{LE1} we
need only to consider the case when $5$ divides either $q-1$ or
$q+1$. By Lemma~\ref{L*}  the Sylow $5$-subgroups of $E_6^-(q)$ are
abelian. Since $5\in\pi(F_4(q))$ there is nothing to prove in this case.

Consider the case $E_6^+ (q)$. Assume first that $5$ divides $q+1$.
Then a Sylow $5$-subgroup of $L$ is a Sylow $5$-subgroup of
$F_4(q)$. Since each element of order $5$ of $F_4(q)$ is real, we have
$5\rightarrow 2\in \E(\Gamma(L))$.

Suppose now that $5$ divides $q-1$. Then the order of a Sylow
$5$-subgroup of $L$ is $5(q-1)_5^6$. By \cite{CaLN} there exists a
torus $T$ of order $\frac{(q- 1)^6}{(q-1,3)}$ with $N_L(T)/T\cong
W(L)$ of order $2^73^45$. Since  $W(L)$ is generated by
involutions, we have that $W(L)$ has no normal $5$-complement and,
in this case, the normalizer of a Sylow $5$-subgroup $S$ of $W(L)$
contains a $2$-element acting non-trivially on $S$. Hence $2\in
\pi(A_5(L))$ and  $5\rightarrow 2\in \E(\Gamma(L))$, as
asserted.\qed

\begin{lem}\label{LE4} $3\rightarrow 2\in \E(\Gamma(L))$ and $\di(\Gamma(L))\leq 4$.
\end{lem}
\noindent {\bf Proof.} Prove first that $3\rightarrow 2\in
\E(\Gamma(L))$. By Lemma~\ref{LE3} we may assume that $r\not=3$.

Suppose that $q\equiv -\epsilon 1\,(\modul 3)$. Then the order of
 a Sylow $3$-subgroup of $L$ is $3^2(q+\epsilon 1)_3^4$.
Since $F_4(q)$ is a subgroup of $L$, we have that every
element of order $3$ of $L$ is real and $3\rightarrow 2\in \E(\Gamma(L))$ by
Lemma~\ref{L4b}.

Now let  $q\equiv \epsilon 1\,(\modul 3)$. Then the order of a
Sylow $3$-subgroup of $L$ is $3^3(q-\epsilon 1)_3^6$. In this case
there is a torus $T$ of $L$ of order $(q-\epsilon 1)^6/(q-\epsilon
1,3)$ and the Sylow $3$-subgroup of $T$  is normalized
by $W(L)$. It is easy to see
that $W(L)$ has self-centralizing Sylow
$3$-subgroups and $2\in \pi(A_3(W))$. By  Lemma~\ref{L0} we
obtain that $3\rightarrow 2\in \E(\Gamma(L))$.

We claim that  $p\rightarrow 2\in \E(\Gamma(L))$ for every prime
$p\in \pi(q^2-1)$. By Lemmas \ref{LE1}, \ref{LE3} and the previous arguments
we may assume that $p\ge 7$. It follows from Lemma~\ref{L*}
that the Sylow $p$-subgroups of $L$ are abelian. Since the elements of
order $p$ in $F_4(q)$ are real, the claim follows  by Lemma~\ref{L4b}.

Consider now  a prime $p\in \pi(\frac{q^6 +\epsilon q^3+1}{(3,q-\epsilon 1)})$.  It follows from \cite{CaLN} and the description
of the centralizers of the elements of order $3$ in L (see lecture by Iwahori  in \cite{BCCI}, and
\cite{Azad}) that there is a Frobenius  subgroup of order $9(q^6+\epsilon q^3+1)$. Hence $p\rightarrow 3$.

Whence,  from Lemmas \ref{LE1}~--~\ref{LE3}, it follows now that $\di(\Gamma(L))\le 4$ and the lemma is proved.\qed

\begin{thm}\label{LE5} Let $G$ be an almost simple group with socle $L=E_6^\pm (q)$.
Then the graph $\Gamma(G)$ is connected and $\di(\Gamma(G))\le 5$.
\end{thm}
\noindent {\bf Proof.} By Lemma~\ref{L2.4} the field automorphisms
of $L$ centralize $W(L)$ and, therefore, for each
semisimple element of order $p$ in $L$, if some edge $p\rightarrow
s$ exists in $\Gamma (L)$, then $p\rightarrow
s\in \E(\Gamma (G))$.
If $r\not=3$, then the edge $r\rightarrow 2$ also remains in
$\Gamma(\Autdiag(L))$.

Suppose that $\varphi $ is a field automorphism of $L$ of order
coprime to $r$. This automorphism preserves all edges
of the form $p\rightarrow s$ if $p\not=r$. On the other
hand,  if $r\not=2$ it is clear that the edge $r\rightarrow 2$
remains in the automorphism group. If $r=2$ we are done by
Lemma~\ref{L3}. Now suppose that $\varphi $ is an
automorphism of $L$ of order coprime to $|L|$. Then, as in
Lemma~\ref{L3.7}, it follows that $d(G)\leq d(L)+1\leq 5$.\qed
\subsection*{Almost simple groups with socle $L=P\Omega _{4l+2}^\epsilon (q)$}
\begin{thm}\label{O} Let $G$ be an almost simple group with  socle $L$ isomorphic to a group $P\Omega _{4l+2}^\epsilon (q)$ ($\epsilon =\pm $, $q=r^t$,
$l\geq 1$). Then the Sylow graph $\Gamma (G)$ is connected and $\di(\Gamma(G))\le 5$.
\end{thm}
\noindent
{\bf Proof.} Define the set $\sigma_1\subseteq \pi(L)$ as follows:\\
If $\varepsilon =+$,
$\sigma_1=\{p\mid  p\mbox{ a prime primitive with respect to the pair }
(q, 2l+1)\}$; \\
if $\epsilon =-$,
$\sigma_1=\{p\mid  p\mbox{ a prime primitive with respect to the pair }
(q, 4l+2)\}$.\\
Set also $\sigma_3:=\{r,2\}$ and $\sigma_2:=(\pi(L)\setminus \sigma_1)\setminus \sigma_3$.
By \cite{St1} we have that $L$ has a subgroup  $M$ isomorphic to $P\Omega _{4l}^\epsilon (q)$ and, moreover,
 $\pi(L)=\pi(M)\cup \sigma_1$.

We prove next that for each odd prime $p\in \sigma_2$ there is an
edge $p\rightarrow 2\in \E(\Gamma (L))$. Clearly this will be also
true for $\Gamma(\Inndiag(L))$ instead of $\Gamma (L)$ by
Lemma~\ref{L3}.

By Lemma~\ref{TV}, since $M$ is a fine group, every semisimple
element of odd prime order in $M$ is real. Let $p\in \sigma_2$.
Since in this case $p\in\pi(M)$ and a Sylow $p$-subgroup of $M$ is a
Sylow $p$-subgroup of $L$, this implies that $p\rightarrow 2\in
E(\Gamma(L))$.

We prove now that for each $p\in \sigma_1$ there exists a prime
$s\in \sigma_2$ such that $p\rightarrow s\in \E(\Gamma(L))$. It is
clear from the Dynkin diagram of $L$ that $A_{2l}\subset
D_{2l+1}$. Hence the group $P\Omega^+ _{4l+2}(q)$ contains the
normalizer of a Singer subgroup of $L_{2l+1}(q)$. This means that
for each prime $p$ primitive with respect to the pair $(q,2l+1)$
and each prime  divisor $s$ of $2l+1$ there is an edge
$p\rightarrow s\in \E(\Gamma(L))$. Clearly, $s\in \sigma_2$. Since
$P\Omega _{4l+2}^-(q)$ contains $U_{2l+1}(q)$, we can deduce as
before that for every prime $p$ primitive with respect to the pair
$(q, 4l+2)$ and $s\in \pi(2l+1)$, there exists an edge
$p\rightarrow s\in \E(\Gamma(L))$.

If $r>2$, the order of the Borel subgroup of $L$ is divisible by
any prime in $\pi(q-1)$ and then $r\rightarrow 2\in
\E(\Gamma(L))$.

It follows from Lemma~\ref{L3} that the graphs $\Gamma(L)$ and
$\Gamma(\Inndiag(L))$ are connected and have diameter at most $4$.

Let $G_0\geq L$ be a maximal normal subgroup of $G$ such that
$\pi(G_0)=\pi(L)$. Clearly, $\Inndiag(L)\leq G_0$. Since the field
automorphisms of $L$ centralize the Weyl group of $L$ (Lemma~\ref{L2.4}), all edges of the form $p\rightarrow s$ for primes in
$\sigma_1\cup \sigma_2$ are also in $\E(\Gamma(G_0))$. Moreover
$r\rightarrow 2\in \E(\Gamma(G_0))$ because every field
automorphism has fixed points on the Cartan subgroup of $L$ of
order at least $\frac{1}{d}(r-1)^{2l}$. Therefore
$\di(\Gamma(G_0))\leq \di(\Gamma (L))\leq 4$. Now we deduce from
Lemma~\ref{L3.7} that $\di(\Gamma(G))\leq 5$.\qed
\subsection*{Almost simple groups with socle $L=L_l^\epsilon(q)$}
\begin{thm}\label{L} Let $G$ be an almost simple group with  socle $L$
isomorphic to $L_l^\epsilon(q)$ ($\epsilon=\pm $, $q= r^t$,
$l>2$). Then the Sylow graph of $G$ is connected and $\di(\Gamma(G))\le 5$.
\end{thm}
\noindent {\bf Proof.} Define the set
$\sigma_1\subseteq \pi(L)$ as follows:\\
If $\epsilon=+$,
 $\sigma_1$ is defined to be the set of primes  primitive with respect to the pair
$(q,j)$ for $j\leq l<2j.$\\
If $\epsilon=-$, $\sigma_1$ is the set of primes $p$ such that
$p\mid q^j-(\epsilon 1)^j\ (j\leq l<2j)$ and $p$ is primitive with
respect to the pair $(q,j)$ when  $j$ is even, and with respect to
the pair  $(q,2j)$ when $j$ is odd.

Set also $\sigma_3:=\{r,2\}$ and $\sigma_2:=(\pi(L)\setminus
\sigma_1)\setminus \sigma_3$.

We prove first that $p\rightarrow 2\in \E(\Gamma(L))$ for any odd prime $p\in \sigma_2\cup \sigma_3$.

If $p=r$, any Sylow $r$-subgroup of $L$  contains  its centralizer
in $\Inndiag(L)$. Hence $r\rightarrow s\in \E(\Gamma(G))$ for
every $s\in\pi(r-1)$. In particular this holds for $s=2$ and
$r>2$.

If $p\in \sigma_2$, then $p$ divides $q^{i}-(\epsilon 1)^i$ for
some $i$, which is minimal with this property. We consider
$l=mi+\lambda$ for some $0\leq \lambda <i$. In this case $m>1$ by
the definition of $\sigma_2$ and there exists a subgroup $R$ of
$K=GL_l^\epsilon(q)$ consisting of all matrices $ \text{diag}
(A_1,\dots, A_m, I)$, where $A_i\in GL_i^\epsilon(q)$ are of order
$f(q,p)=(q^i-(\epsilon 1)^i)_p$ and $I$ is the identity
$\lambda\times \lambda$-matrix.

It follows from arithmetical arguments that $|K|_p=(q^i-(\epsilon
1)^i)_p^m(m!)_p$. The subgroup $R$ is a direct product of $m$
copies of groups of order $f(q,p)$. Obviously the subgroup $S_m$
permuting blocks of the matrices of size $i$ is contained in
$SL_l^\epsilon(q)$ and normalizes $R$. Moreover, this subgroup
contains a Sylow $p$-subgroup of $K$. Now we can deduce from
Corollary~\ref{C2} that $2\in \pi(A_p(L))$ and so $p\rightarrow
2\in \E(\Gamma (L))$.

We claim that for every $p\in \sigma_1$  there exists an
edge $p\rightarrow s$ for some  $s\in \sigma_2\cup \sigma_3$.
Let $p\in \sigma_1$. If $p\mid
q^j-(\epsilon 1)^j$ with $2j>l$, then a Sylow $p$-subgroup of
$L$ is cyclic and isomorphic to a subgroup of a Singer subgroup
of the group $L_j^\epsilon(q)$. Since $2j>l$, this implies, from
properties of  Singer subgroups, that $|L|_p=(q^j-(\epsilon 1)^j)_p$. Moreover, in this case,
for every prime divisor $s$ of $j$ there exists an element of
order $s$ that normalizes, but not centralizes, a Sylow $p$-subgroup of $L$. Hence there
exists an edge $p\rightarrow s$. We see next that $s\in \sigma_2\cup \sigma_3$. This is clear if either $j\equiv 0\,(\modul r)$ or
$j\equiv 0\,(\modul 2)$. Hence we may assume that $j$ is odd. If $j$ is not a prime,
since $s$ divides $j$, we have $s\mid q^{s-1}-1$, $s\leq
\frac{j}{2}$ and hence $2s\leq j\leq l$. This forces $2(s-1)<l$
and $s\in \sigma_2$.
Suppose that $j=s$. As above $s\mid q^{s-1}-1$. If $s$ is
primitive with respect to the pair $(q,s-1)$,
then $s\notin \sigma_1$ and so $s\in \sigma_2\cup \sigma_3$.
Otherwise, $s\mid q^i-1$ for some $i$ dividing
$s-1$. Hence $2i<l$ and $s\in \sigma_2$, which proves the claim.

Whence it follows that in fact the group
$PGL_l^\epsilon(q)$ has Sylow graph with diameter at most $4$.
We can argue now like in the proof of Theorem~\ref{O}, for $L=P\Omega_{4l+2}^\pm(q)$, to obtain finally the desired result.\qed

\medskip
The Main Theorem follows now from Corollary~\ref{C2} and
Theorems~\ref{P1}, \ref{S}, \ref{fs}, \ref{LE5}, \ref{O}, \ref{L}.

\section{An application: Sylow normalizers}\label{Sylnorm}

In this section we apply the results on the Sylow graph of a group
to obtain some information about the structure of groups.
\begin{de} (\cite[Definition]{DGP}) A \emph{covering-formation} $\mathcal F= LF(f)$ is a non-empty saturated formation
 locally defined by a formation
function $f$ given by:
\[f(p)=\mathcal E_{\pi(p)}\mbox{ if } p\in \pi,\quad
f(q)=\emptyset \mbox{ if } q\not\in \pi,\] $\pi:=\Char(\mathcal
F)$, the characteristic of $\mathcal F$, where for each prime
$p\in \pi$, $\pi(p)$ is a set of primes satisfying the following
conditions:

(i) $p\in \pi(p)\subseteq \pi$,

(ii)  $q\in \pi(p)\implies~  p\in \pi(q)$, \\
and where $\mathcal E_{\pi(p)}$ is the class of
all finite $\pi(p)$-groups.\\
\emph{(As mentioned in \cite[Remark 1(b)]{DGP} the condition $\pi(p)\subseteq \pi$
is no loss of generality.)}

\medskip
Let $\mathcal S$ denote the class of all finite soluble groups and
let $\mathcal S_{\pi(p)}$ be the class of all finite soluble
$\pi(p)$-groups. The saturated formation $\mathcal F\cap \mathcal
S=LF(g)$ such that
$$g(p)=\mathcal S_{\pi(p)}\mbox{ if } p\in \pi,\quad
g(q)=\emptyset \mbox{ if } q\not\in \pi,$$  is called a
\emph{covering-formation of soluble groups}.
\end{de}

With the previous notation, if $q\in \pi(p)$, we write $p\leftrightarrow  q$. (This defines a reflexive and symmetric relation on $\pi$.)
We associate to the covering-formation $\mathcal F$ another family of sets of primes defined as follows:
\[\Sigma _\mathcal F:=\{\sigma \subseteq \pi \mid \mbox{ if } p,q\in \sigma \mbox{ and } p\not=q,\mbox{ then } p\not\leftrightarrow q\}.\]

We say that two primes $r$ and $p$ of $\pi$ are \emph{connected} if either $r=p$ or there is a sequence of primes
$r=r_1,\ldots,\,r_n=p$ in $\pi$
such that $r_i\leftrightarrow  r_{i+1}$ for any $1\leq i\leq n-1$.

\medskip
\emph{In the sequel $\mathcal F$  always denotes a covering-formation
as defined above.} We introduce some more notation.

For a group $G$ and a prime $p$, let $G_p$  denote a Sylow
$p$-subgroup of $G$. Recall that $\pi(G)$ denotes the set of all
primes dividing the order of $G$. Then, for a class $\mathcal X$
of groups, the class map {\scshape n} is defined as follows:
\[\mbox{{\scshape n}}\mathcal X=(G\mid N_G(G_p)\in \mathcal X,\mbox{ for every prime }p\in \pi (G)).\]
For notation and results about classes of groups and closure operations we refer to Doerk and Hawkes book~\cite{DH}.

\smallskip
We are interested in the study of classes of groups $\mathcal X$ satisfying
$\mbox{{\scshape n}}\mathcal X= \mathcal X$, that is, a group $G$ belongs to such a class $\mathcal X$ if and only if its Sylow normalizers belong to the class.

 \smallskip
 In \cite{DGP} the covering-formations of soluble groups are characterized as
the subgroup-closed saturated formations $\mathcal X$ satisfying that $\mbox{\scshape n}\mathcal X=\mathcal X$
in the universe of all finite soluble groups. Also they are exactly the classes of soluble groups with nilpotent Hall subgroups for prefixed sets of primes.
For any set of primes $\tau$ we set $\mathbf E_\tau^{\n}$ for the class of all groups with nilpotent Hall $\tau$-subgroups. To be more precise we gather the following results.
\begin{prop}\label{s1} \emph{\cite[Proposition 2]{DGP}}
Let $\mathcal G = \mathcal F\cap \mathcal S$ be a covering-formation of soluble groups as defined above. Then:
\begin{enumerate}
\item $\mbox{\scshape s}\,{\cal G} = {\cal G}\mbox{ and {\scshape n}}\,{\cal G}\cap \mathcal S={\cal G};$
\item $\mathcal G=\bigcap_{\sigma\in \Sigma_\mathcal F,|\sigma|=2}\mathbf E_\sigma^{\n}\cap \mathcal S=\bigcap_{\sigma\in
\Sigma_\mathcal F}\mathbf E_\sigma^{\n}\cap \mathcal S.$
\end{enumerate}

\end{prop}
\begin{thm}\label{s2} \emph{\cite[Theorem]{DGP}}
Let $\cal H$ be a
 subgroup-closed saturated
formation of soluble groups. Then the following statements are equivalent:

(i) For any  soluble group, its Sylow normalizers belong to $\mathcal H$ if and only if the group belongs to $\mathcal H$;

(ii) $\mathcal H$ is a covering-formation of soluble groups.
\end{thm}

The next technical result will be useful to clarify the
constructions  in the examples below.
\begin{lem}\label{lem0} Let $\{\tau_i\mid i\in I\}$ be a partition of a set of primes $\tau$ and let
$\mathcal D=LF(d)$ be a saturated formation of
characteristic $\tau$, locally defined by a formation function $d$ satisfying  that $\mathcal S_pd(p)\subseteq \mathcal E_{\tau_i}$
for every $i\in I$ and every prime $p\in \tau_i\subseteq \tau$. Then
$$\mathcal D=\times_{i\in I}\mathcal D_i  :=(G\in \mathcal E_\tau \mid
G=\times_{i\in I}G_{\tau_i},\ G_{\tau_i}\in \mathcal D_i)$$
where for every $i\in I$, $\mathcal D_i:=\mathcal D\cap \mathcal E_{\tau_i}=LF(d_i)$ is the saturated formation of characteristic $\tau _i$,
locally defined by the formation function $d_i$
given by $$d_i(p)=d(p)\mbox{ if } p\in \tau_i\ ,\quad d_i(q)=\emptyset \mbox{ if } q\not\in \tau_i.$$

\end{lem}
\noindent
{\bf Proof.} The result is easily proven taking into account that
\begin{gather*}
\mathcal D= \mathcal E_\tau\cap (\cap_{p\in \tau}\mathcal E_{p'}\mathcal S_p d(p))\subseteq
\mathcal E_\tau\cap (\cap_{i\in I}\mathcal E_{\tau_i'}\mathcal E_{\tau_i})=
(G\in \mathcal E_\tau\mid O^{\tau_i}(G)\in \mathcal E_{\tau_i'},\
i\in I)=\\=(G\in \mathcal E_\tau \mid G=\times_{i\in I}G_{\tau_i},\ G_{\tau_i}\in \mathcal E_{\tau_i}).\qed
\end{gather*}
\begin{rem*} \cite[Proposition IV(3.8)]{DH} We recall that if $\mathcal D=LF(D)=LF(d)$ is a saturated formation with canonical local definition $D$
and where $d$ is any other local
definition of $\mathcal D$,  then
 $D(p)=\mathcal S_pd(p)\cap \mathcal D$
for any prime $p$.

 Hence the hypothesis on the local definition $d$ in Lemma~\ref{lem0} can be equivalently replaced by the assumption that the canonical local definition $D$ of $\mathcal D$ satisfies that $D(p)\subseteq \mathcal E_{\tau_i}$
for every $i\in I$ and every prime $p\in \tau_i\subseteq \tau$.
\end{rem*}
\begin{ex} Let  $\{\pi_i\mid i\in I\}$ be a partition of the set of primes $\pi$ and
$$\times_{i\in I}\mathcal E_{\pi_i}  :=(G\in \mathcal E_\pi \mid G=\times_{i\in I}G_{\pi_i},\ G_{\pi_i}\in \mathcal E_{\pi_i}).$$
As  follows  from Lemma~\ref{lem0}, this class is  a
covering-formation $\mathcal F=\times_{i\in I}\mathcal E_{\pi_i}$
such  that the sets of primes $\pi(p)$ ($p\in \pi$) form the given
partition of $\pi$.
 \end{ex}

\medskip
Example 1 is a particular type of the following construction:

\begin{ex}\label{ex2} Let $\mathcal H$ be a covering-formation of (soluble) groups of characteristic $\rho $
($2\not\in \rho  $), $\pi=\rho\cup \tau$ and $\rho\cap
\tau=\emptyset$. Let
\[\mathcal H\times \mathcal E_\tau :=(G\in \mathcal E_\pi\mid G=A\times B,\ A\in \mathcal H,\ B\in \mathcal E_\tau ).\]
Then $\mathcal F=\mathcal H\times \mathcal E_\tau$ is a covering-formation  such that if $2, p\in
\pi$ and $p$ is connected with $2$, then $\pi(p)=\pi(2)=\tau$.
Moreover, $\mathcal F=\mathcal H\times \mathcal E_\tau
=\bigcap_{\sigma\in \Sigma_\mathcal F}\mathbf E_\sigma^{\n}.$
\end{ex}
\noindent
{\bf Proof.} The first part follows from Lemma~\ref{lem0}. Then from Proposition~\ref{s1} we can deduce that
$\mathcal F\subseteq \mathcal L:=\bigcap_{\sigma\in \Sigma_\mathcal F}\mathbf E_\sigma^{\n}.$

Assume now that $\mathcal F$ is properly contained in $\mathcal L$ and let $G\in \mathcal L\setminus \mathcal F$ of minimal order. Since
$\mathcal L$ is {\scshape q}-closed and $\mathcal F$ is a saturated formation, $G$ is a primitive group of type either $1$ or $2$, with a unique minimal
normal subgroup, say $N$, and $G/N\in \mathcal F$.

Assume that  $G$ is primitive of type $1$ and $N$ is a $p$-group for a prime $p$. Since $G\in \mathcal L$, it follows that
$G/N=G/C_G(N)\in \mathcal E_{\pi(p)}=f(p)$ and so $G\in \mathcal F$, a contradiction.

Therefore $G$ is a primitive group of type $2$ and, in particular,
$2\in \pi(N)$. Note that any normal subgroup of $G$, in particular
$N$, belongs to $\mathcal L$.

It is known that no simple group $E$ has nilpotent Hall $\nu
$-subgroups with $\nu $ a set of primes such that $|\nu \cap
\pi(E)|>1$ and $2\in \nu $  (see \cite[Proposition 1]{AC}). Hence
any non-abelian simple group in $\mathcal L$ belongs to $\mathcal
E_{\pi(2)}$; consequently, $N\in \mathcal E_{\pi(2)}\subseteq
\mathcal F$ and $N<G$.

On the other hand, the choice of $G$ implies that $M\in \mathcal F\subseteq \mathcal E_{2'}\mathcal E_{\pi(2)}$ for
a maximal normal subgroup $M$ of $G$.
Consequently $M\in \mathcal E_{\pi(2)}$ because $N\leq M$ and then $O_{2'}(M)=1$.

Since $G\not\in \mathcal E_{\pi(2)}$, there is a prime $p\in \pi(G)\setminus \pi(2)$ and $|G/M|=p$
as $G/M\in \mathcal F=\mathcal H\times \mathcal E_\tau$. Since $G\in \mathcal L$ it follows that $G$ has a normal
Sylow $p$-subgroup $G_p$ and $G=M\times G_p\in \mathcal F$, the final contradiction.\qed

\begin{ex}
Let $\mathcal F$ be a covering-formation such that if $2,p\in \pi$
and  $p\in \{3,5\}$, then $p$ is not connected with $2$.
(Equivalently, $\pi=\rho \cup \tau $ with $\rho \cap
\tau=\emptyset$, $2\not\in \rho$ and $3,5\not\in \tau$). Then
$\mathcal F\subseteq \mathcal S$ is a covering-formation of
soluble groups. More precisely,
$$\mathcal F=\mathcal H\times \mathcal J:=(G\in \mathcal E_\pi\mid G=H\times J,\ H\in \mathcal H,\ J\in \mathcal J)$$
where $\mathcal H$ and $\mathcal J$ are  covering-formations of characteristics $\tau$ and $\rho$, respectively.
Moreover, $\mathcal F=\mathcal H\times \mathcal J=\bigcap_{\sigma\in \Sigma_\mathcal F}\mathbf E_\sigma^{\n}.$\\
\emph{(For instance the following classes of groups appear as particular types of this construction:
\[\mathcal H\times \mathcal E_{\{3,5\}}:=(G\in \mathcal E_\pi\mid  G=A\times B,\ A\in \mathcal H,\ B\in \mathcal E_{\{3,5\}}),\mbox{ or also}\]
\[\mathcal H\times \mathcal E_3\times \mathcal E_5:=(G\in \mathcal E_\pi\mid  G=A\times
B\times C,\ A\in \mathcal H,\ B\in \mathcal E_{3},\ C\in \mathcal
E_5),\] where (in both cases) $3,5\in \pi$ and $\mathcal H$ is a
covering-formation of soluble groups of characteristic
$\pi\setminus \{3,5\}$.)}

\end{ex}
\noindent {\bf Proof.} The restrictions on the sets of primes
imply that $f(p)=\mathcal E_{\pi(p)}\subseteq \mathcal S$ for any
$p\in \pi$. Hence $\mathcal F\subseteq \mathcal S$. See
Lemma~\ref{lem0} for the description of $\mathcal F$.

On the other hand,  by Proposition~\ref{s1} we have $\mathcal
F=\mathcal L\cap \mathcal S$ for $\mathcal L:=\bigcap_{\sigma\in
\Sigma_\mathcal F}\mathbf E_\sigma^{\n}$. Let us consider the
covering-formation $\mathcal C:=\mathcal E_\tau\times \mathcal J$,
which is of type given in  Example~\ref{ex2}. Then  $\mathcal C=
\bigcap_{\sigma\in \Sigma_\mathcal C}\mathbf E_\sigma^{\n}$. Note
that $\Sigma_\mathcal C\subseteq \Sigma_\mathcal F$, which implies
that $\mathcal L=\bigcap_{\sigma\in \Sigma_\mathcal F}\mathbf
E_\sigma^{\n}\subseteq \bigcap_{\sigma\in \Sigma_\mathcal
C}\mathbf E_\sigma^{\n}= \mathcal C\subseteq \mathcal S.$
 Consequently $\mathcal F=\mathcal L$ and we are done.\qed

\begin{rem}\label{rem1} (1) For covering-formations, in the universe of all finite groups, we know that $\mathcal F\subseteq \mbox{{\scshape n}}\mathcal F\subseteq \mathcal E_\pi$ since  $\mathcal F$ is
subgroup-closed, but in general $\mathcal F\not= \mbox{{\scshape n}}\mathcal F$. (See \cite[Remark 1(c)]{DGP}.)
\smallskip

(2) A subgroup-closed saturated formation $\mathcal X$ with
$\mbox{{\scshape n}}\mathcal X=\mathcal X$ need not be a
covering-formation. (See \cite{DDG}.)
\smallskip

(3) As pointed out also in \cite{DGP}, it follows from the
Feit-Thompson theorem and Theorem~\ref{s2} that if $2\notin \pi$,
then $\mathcal F= \mbox{{\scshape n}}\mathcal F$. The result
$\mathcal F= \mbox{{\scshape n}}\mathcal F$ can be deduced in the
same way under the hypothesis that $3,5\not\in \pi$.
\end{rem}

\begin{prop}\label{mincon} Let $G$ be a group of minimal order in $\mbox{{\scshape n}}\mathcal F\setminus \mathcal F$.
Then $G$ is a primitive group of type $2$ and any pair of primes
in $\pi(G)$ is connected.

\end{prop}
\noindent {\bf Proof.} Note that $\mbox{{\scshape n}}\mathcal F$
is {\scshape q}-closed because $N_{G/N}(G_pN/N)=N_G(G_p)N/N$ for
any $G_p\in  \Syl_p(G)$ and  $ N\unlhd G$. Since $\mathcal
F$ is a saturated formation we can deduce that $G$ is a primitive
group of type either $1$ or $2$, i.e., with a unique minimal
normal subgroup, say $N$.

Assume that $G$ is primitive of type 1. Then $N=C_G(N)$ is a $p$-group, for some prime $p$, and $G=NM$ with $M$ a maximal
subgroup of $G$, $N\cap M=1$.

We claim that $G/N\in \mathcal E_{\pi(p)}=f(p)$.

If $p\not\in \pi(G/N)$, then $N\in \mbox{Syl}_p(G)$ and
$G=N_G(N)\in \mathcal F$, a contradiction.

Hence $p\in \pi(G/N)$. Assume that there exists $q\in \pi(G/N)$,
$q\not\in \pi(p)$. Observe that $M\cong G/N\in \mathcal F\subseteq
\mathcal E_{p'}\mathcal E_{\pi(p)}$. In particular we can consider
$1\not=M_q\in \mbox{Syl}_q(M)$ and $M_q\leq O_{p'}(M)$. By Lemma~\ref{L0} it follows that $M=O_{p'}(M)N_M(M_q)$. Then
there exists $1\not=M_p\in \mbox{Syl}_p(M)$ with $M_p\leq
N_M(M_q)\in \mathcal F$. Since $p\not\in \pi(q)$ it follows that
$[M_p,M_q]=1$. Consequently $M_q\leq N_G(NM_p)\in \mathcal F$ as
$NM_p\in \mbox{Syl}_p(G)$ and $G\in \mbox{{\scshape n}}\mathcal
F$. Again since $q\not\in \pi(p)$ it follows that $M_q\leq
C_G(N)=N$, a contradiction which proves the claim.

Since $G/N\in \mathcal F$ and $G/N=G/C_G(N)\in \mathcal E_{\pi(p)}=f(p)$ it
follows that $G\in \mathcal F$, a contradiction.

Therefore $G$ is a primitive group of type $2$. Hence
$N=N_1\times\ldots \times N_r$, where $N_i\cong N_j$ are non-abelian
simple groups for any $i,j\in \{1,\ldots,r\}$. Set
$R=\cap_{i=1}^{r} N_G(N_i)\unlhd G$.

Assume first that $R<G$. We claim that if $p\in \pi (N)$ and $q\in
\pi(G/R)$, then $p\leftrightarrow  q$. This will imply the result.

Since $G/R\in \mathcal F$ we have, in particular, that $G/R\in
\mathcal E_{q'}\mathcal E_{\pi(q)}$. If $p\not\in \pi(q)$ there
exists $G_p\in \mbox{Syl}_p(G)$ such that $G_p\leq T\unlhd G$ with
$T/R=O_{q'}(G/R)$. By Lemma~\ref{L0}, $G=TN_G(G_p)$ and there
exists a subgroup $T_qQ\in \mbox{Syl}_q(G)$ for some $T_q\in
\mbox{Syl}_q(T)$ and some $Q\in \mbox{Syl}_q(N_G(G_p))$. Since
$N_G(G_p)\in \mathcal F$ and $q\not\in \pi(p)$ it follows that
$[G_p,Q]=1$ and thus $[N_p,Q]=1$, where $1\not=N_p=G_p\cap N\in
\mbox{Syl}_p(N)$. But $R<G$. Hence  $r>1$ and
$N_p=(N_1)_p\times\ldots \times(N_r)_p$ with $(N_i)_p\in
\mbox{Syl}_p(N_i)$ for every $i=1,\ldots,r$. Since $G$ acts
transitively via conjugation on the set $\{N_1,\ldots,N_r\}$, we
deduce that $Q\leq R$ and therefore $T_qQ\leq R$, a contradiction.

Assume now that $G=R$. Then $r=1$ and $G$ is an almost simple group.
We claim that $\pi(A_p(G))\subseteq \pi(p)$ for every prime $p\in \pi(G)$.
Since $N_G(G_p)\in \mathcal F$ for $G_p\in \mbox{Syl}_p(G)$, if $q\in \pi(N_G(G_p))$ and $q\not\in \pi(p)$
we have that
$[Q,G_p]=1$ for any $Q\in \mbox{Syl}_q(N_G(G_p))$. Therefore $q\not\in \pi(A_p(G))$ and the claim is proved.
The result can be deduced now from the Main Theorem.\qed

\medskip
This proposition can be useful to derive conditions assuring that $\mbox{{\scshape n}}\mathcal F=\mathcal F$. In particular we deduce next that
for any of the classes of groups in Examples 1, 2 and 3, a group belongs to the class if
and only if its Sylow normalizers belong to the class.

\begin{cor}\label{cor0} For any of the classes of groups in Example 1, a group belongs to the class if
and only if its Sylow normalizers belong to the class.
\end{cor}

\noindent
{\bf Proof.} Let $\mathcal F$ be a covering-formation as in Example 1. Since $\mathcal F$ is subgroup-closed, we need only to prove that
 $\mbox{{\scshape n}}\mathcal F\subseteq \mathcal F$. Assume that this is
wrong and let $G$ be a group of minimal order in $\mbox{{\scshape n}}\mathcal F\setminus
\mathcal F$. Then Proposition~\ref{mincon}
implies that $G$ is a primitive group of type $2$ and each pair of
primes $p,q\in \pi(G)$ is connected. In particular, every prime
$p\in \pi (G)\setminus\{2\}$ is connected with $2\in \pi (G)$ and
consequently $\pi(p)=\pi(2)$ and $G\in \mathcal
E_{\pi(2)}\subseteq \mathcal F$, a contradiction.\qed
\begin{rem*} Observe that the proof of Corollary~\ref{cor0} can be adapted
to derive the corresponding results for the covering-formations in
Examples 2 and 3. We show below a more direct approach to these
results based on Corollary~\ref{cor0}.
\end{rem*}
\begin{lem}\label{lem1}
Let $\{\tau_i\mid i\in I\}$ be a partition of a set of primes $\tau$ and let $\mathcal X$ and $1\in \mathcal X_i\subseteq \mathcal E_{\tau_i}$
 for every $i\in I$ be classes of groups such that
$$\mathcal X:=\times_{i\in I}\mathcal X_i  :=(G\in \mathcal E_\tau \mid
G=\times_{i\in I}G_{\tau_i},\ G_{\tau_i}\in \mathcal X_i).$$
 Then {\scshape n}$\mathcal X=\mathcal X$ if and only if {\scshape n}$\mathcal X_i=\mathcal X_i$ for all $i\in I$.
\end{lem}
\noindent
{\bf Proof.} This follows by an straightforward proof and Corollary~\ref{cor0}.\qed
\begin{cor}\label{cor2-3} For any of the classes of groups in Examples 2 and 3, a group belongs to the class if
and only if its Sylow normalizers belong to the class.
\end{cor}
\noindent
{\bf Proof.} It follows by Lemma~\ref{lem1} and Remark~\ref{rem1} (3).\qed

\begin{rem*} As pointed out in Remark~\ref{rem1} (2), there is an example
in \cite{DDG} of a  subgroup-closed saturated formation $\mathcal
X$ with the property $\mbox{{\scshape n}}\mathcal X= \mathcal X$
which is not a covering-formation. This example depends on the
primes in the characteristic of $\mathcal X$ and on the non-abelian
simple groups involving these primes as divisors of their orders.
Here $\sigma:=\Char (\mathcal X)=\{2,5,7,13,p\}$, where $p\not=2,3,5,7,13$,
and it is proved that the Suzuki group $\Sz(2^3)$ is the only
non-abelian simple $\sigma$-group.

We may wonder whether there is an example of such a formation with full characteristic. Lemma~\ref{lem1}
provides an easy way to
construct an example as aimed from the saturated formation $\mathcal X$ in \cite{DDG}. We need just to consider now the subgroup-closed
saturated formation $\mathcal H:=\mathcal X\times \mathcal E_{\sigma'}$ with $\sigma=\Char (\mathcal X)$.
\end{rem*}

\end{document}